%% Paper intended for special issue of "Publication Inst. Math."
%% dedicated to the centenary of Jovan Karamata (1902-1967).
%% This is an Amstex file, version of 11.02.2001, final.
%% ATTENTION: if there are font problems, replace `leqslant' and
%% `geqslant' by `le' and `ge'.
%% In print in "Publs. Inst. Math." 71(85)(2002)

\input amstex.tex
\documentstyle{amsppt}
%\magnification=1200
%\hsize=12.5cm
%\vsize=19cm
%\hoffset=1cm
%\voffset=2cm
\def\DJ{\leavevmode\setbox0=\hbox{D}\kern0pt\rlap
{\kern.04em\raise.188\ht0\hbox{-}}D}
\footline={\hss{\vbox to 2cm{\vfil\hbox{\rm\folio}}}\hss}
\nopagenumbers
\font\ff=cmr8

\baselineskip=13pt

\def\a{\alpha}\def\b{\beta}
\def\d{{\,\roman d}}
\def\e{\varepsilon}

\def\={\;=\;}

\def\D{\Delta}

\def\z{\zeta}

\def\e{\varepsilon}

\def\e{\varepsilon}
\def\l{\lambda}

\font\seveneufm=eufm7
\font\fiveeufm=eufm5
\newfam\eufmfam
\textfont\eufmfam=\teneufm
\scriptfont\eufmfam=\seveneufm
\scriptscriptfont\eufmfam=\fiveeufm
\def\mathfrak#1{{\fam\eufmfam\relax#1}}

\font\tenmsb=msbm10
\font\sevenmsb=msbm7
\font\fivemsb=msbm5
\newfam\msbfam
\textfont\msbfam=\tenmsb
\scriptfont\msbfam=\sevenmsb
\scriptscriptfont\msbfam=\fivemsb
\def\Bbb#1{{\fam\msbfam #1}}

\def \RR {\Bbb R}

\def\rightheadline{{\hfil{\ff
Additive functions and regular variation}\hfil\tenrm\folio}}

\def\leftheadline{{\tenrm\folio\hfil{\ff
A. Ivi\'c }\hfil}}
\def\emptyheadline{\hfil}
\headline{\ifnum\pageno=1 \emptyheadline\else
\ifodd\pageno \rightheadline \else \leftheadline\fi\fi}

\topmatter
\title ON SUMMATORY FUNCTIONS OF ADDITIVE FUNCTIONS AND
REGULAR VARIATION
\endtitle
\author   Aleksandar Ivi\'c \endauthor
\address{ \bigskip
Aleksandar Ivi\'c, Katedra Matematike RGF-a
Universiteta u Beogradu, \DJ u\v sina 7, 11000 Beograd,
Serbia (Yugoslavia). \bigskip}
\endaddress
\keywords Regularly and slowly varying functions, arithmetic sums,
additive functions, Abelian and Tauberian theorems
\endkeywords
\subjclass Primary 11N37,  secondary 26A12 \endsubjclass
\email {\tt aivic\@matf.bg.ac.yu,
aivic\@rgf.bg.ac.yu} \endemail
\abstract
{An overview of results and problems concerning the
asymptotic behaviour for summatory functions
of a certain class of additive functions is given.
The class of functions in question involves regular variation. Some
new Abelian and Tauberian results for additive functions of the
form $F(n) = \sum_{p^\a||n}\a h(p)$ are obtained.}
\endabstract
\endtopmatter

\noindent

\def\a{\alpha}
\def\e{\varepsilon}

\def\={\;=\;}

\font\teneufm=eufm10
\font\seveneufm=eufm7
\font\fiveeufm=eufm5
\newfam\eufmfam
\textfont\eufmfam=\teneufm
\scriptfont\eufmfam=\seveneufm
\scriptscriptfont\eufmfam=\fiveeufm
\def\mathfrak#1{{\fam\eufmfam\relax#1}}

\font\tenmsb=msbm10
\font\sevenmsb=msbm7
\font\fivemsb=msbm5
\newfam\msbfam
      \textfont\msbfam=\tenmsb
      \scriptfont\msbfam=\sevenmsb
      \scriptscriptfont\msbfam=\fivemsb
\def\Bbb#1{{\fam\msbfam #1}}

\def \RR {\Bbb R}

 \def\D{\Delta}  \def\e{\varepsilon}
 
\def\DJ{\leavevmode\setbox0=\hbox{D}\kern0pt\rlap
 {\kern.04em\raise.188\ht0\hbox{-}}D}

\heading 1. Introduction
\endheading

The purpose of this paper is to provide an overview of results and
problems concerning the
asymptotic behaviour of summatory functions
for a certain class of additive functions.
An arithmetic function $g(n)$ is {\it additive} if $g(mn) = g(m) + g(n)$
whenever $(m,n) = 1$. Many additive functions can be represented in the form
$$
f(n) = \sum_{p|n}h(p)\quad\text{or}\quad
F(n) = \sum_{p^\a||n}\a h(p),\quad h(p) = p^\rho L(p),\leqno(1.1)
$$
where $\rho\in\RR$ and $p^\a||n$ means ($p$ denotes primes) that
$p^\a$ divides $n$, but $p^{\a+1}$ does not. The function $h(x)$  appearing
in (1.1)  is a {\it regularly varying} function. It is a positive, continuous
function for $x \geqslant x_0 \,(>0)$, for which
there exists $\rho \in\RR$ (called the index of $h$) such that
$$
\lim_{x\rightarrow \infty} \frac{h(cx)}{h(x)}=c^{\rho},
\quad \text{for all }\; c>0. \leqno(1.2)
$$
We shall denote the set of all regularly varying functions
by $\Re$. We shall also denote by ${\Cal L}$ the set of
slowly  {\it slowly varying} (or {\it slowly oscillating})
functions, namely those functions in $\Re$
for which the index $\rho=0$. It is easy to show that
if $h\in \Re$, then there exists $L\in {\Cal L}$
such that $h(x)=x^{\rho}L(x)$, with $\rho$ being the index of $h$.

Slowly varying functions arise naturally in many branches of number
theory, especially in the theory of arithmetic functions.   Namely
if $a(n)$ is an arithmetic function, then the summatory function
$A(x) := \sum_{n\leqslant x}a(n)$ is usually written in the form
$$
A(x) = M(x) + E(x),\leqno (1.3)
$$
where $M(x)$ is the {\it main term} and $E(x)$ the {\it error term} in
the asymptotic formula (1.3), meaning
$$
\lim_{x\to\infty}\,{E(x)\over M(x)} \= 0.                 \leqno(1.4)
$$
The main term $M(x)$ is usually represented in terms of elementary
functions such as powers, exponentials, logarithms etc.
The error term function $E(x)$ is usually not smooth, but is bounded in
terms of smooth functions, and one often looks for estimates of the form
$$
E(x) \= O(\D(x)),      \leqno(1.5)
$$
where $\D(x)$ is a regularly varying function.

For a comprehensive account of regularly varying functions
the reader is referred to the
monographs of Bingham et al. [1] and E. Seneta [21]. By a fundamental
result of J. Karamata [20], who founded the theory of regular variation,
the limit in (1.2) is uniform for $0 < a \leqslant c \leqslant b < \infty$
and any $0 < a < b$. This is known as the {\it uniform convergence
theorem}. It is used to show that any slowly
varying function $L(x)$ is necessarily of the form
$$
L(x) = A(x)\exp\left(\int_{x_0}^x\eta(t){\d t\over t}\right),\;
\lim_{x\to\infty}A(x) = A > 0,\;   \lim_{x\to\infty}\eta(x) = 0.\leqno(1.6)
$$
Usually the function $\D(x)$ in (1.5) is assumed to be of the form
$\D(x) = x^\rho L(x)$, where $L(x)$ is of the form (1.6) with
$A(x) \equiv A \,(>0)$. This is convenient, since we are interested in
the asymptotic behaviour of $\D(x)$ and  $\lim_{x\to\infty}A(x) = A$
has to hold in (1.6).

Special cases of (1.1) (with $\rho = 0, L(x) \equiv 1$) are
the well-known functions
$$
\omega(n) = \sum_{p|n}1,\qquad \Omega(n) = \sum_{p^\a||n}\a,  \leqno(1.7)
$$
which represent the number of distinct prime factors of $n$ and the
total number of prime factors of $n$, respectively.
Another pair of additive functions, which in distinction with the
``small additive" functions of (1.7)  are called ``large additive functions",
is (see e.g., [10], [12] and [19])
$$
\b(n) = \sum_{p|n}p,\qquad B(n) = \sum_{p^\a||n}\a p.  \leqno(1.8)
$$
This pair is the special case $\rho = 1,\,L(x) \equiv1$ of (1.1).

Given an additive arithmetic function $g = f$ or $g = F$, defined by (1.1),
one can ask for properties of the summatory function
$$
G(x) \;\=\; \sum_{n\leqslant x}g(n)
$$
from the properties of $g(n)$. Such types of results are
commonly called {\it Abelian theorems}. Conversely, if from the properties
of $G(x)$ we can deduce properties of $g(n)$ itself, then
these kinds of results are called {\it Tauberian theorems}.
Abelian and Tauberian theorems appear in various contexts in many
branches of analysis and, in general, it is Tauberian theorems
which are more difficult to prove than their corresponding
Abelian counterparts. Tauberian theorems commonly necessitate an
additional so-called {\it Tauberian condition}, which usually
involves monotonicity of some kind or non-negativity.

\heading
2. Abelian  theorems for sums of additive functions
\endheading

If $f(n)$ is given by (1.1) with $L\in{\Cal L}$ and $\rho > 0$, then
we have
$$
\sum_{n\leqslant x}f(n) = \left({\z(\rho+1)\over\rho+1}+o(1)\right)
{x^{\rho+1}L(x)\over\log x}\quad(x\to\infty),\leqno(2.1)
$$
where $\z(s)$ is the Riemann zeta-function (see [15]). This Abelian
theorem was proved in [6], with the remark that the very generality of
$L(x)$ hinders a result sharper than (2.1), namely to have
an $O$--term in place of $o(1)$. The asymptotic formula
(2.1) cannot hold for $\rho = 0$, since $\z(s)$ has a pole at $s=1$.
In this case the sum in question equals ($[y]$ is the greatest integer
not exceeding $y$)
$$
\int_2^x{L(t)\over\log t}\left[{x\over t}\right]\d t\leqno(2.2)
$$
plus an error term, which depends on the prime number theorem (see [15,
Chapter 12]). Bingham--Inoue [3] considered the case $\rho=0$ in (2.1) and
obtained that
$$
\sum_{n\leqslant x}\sum_{p|n}L(p) \;\sim\; xg(x)
\qquad(x\to\infty),\leqno(2.3)
$$
where
$$
g(x) \;:=\; \int_1^{x/2}{[u]\over u^2}{\tilde\ell}\left({x\over u}\right)\d u,
\qquad {\tilde \ell}(x) \;:=\; {L(x)\over\log x}.
\leqno(2.4)
$$
Moreover, $g(x) \in \Pi_{\tilde \ell}$ with ${\tilde \ell}$--index 1,
where the {\it de Haan class}  $\Pi_\ell$, for given
$\ell\in{\Cal L}$, consists of measurable functions $g$ satisfying
$$
\lim_{x\to\infty}\,{g(\lambda x)-g(x)\over \ell(x)} \= c\log\lambda
\qquad(\forall \lambda > 0,\,x\to\infty),\leqno(2.5)
$$
where $c$ is called the $\ell$-{\it index} of $\Pi_\ell $.
We shall give now a direct proof of the fact that
$g(x) \in \Pi_{\tilde \ell}$ with ${\tilde \ell}$--index 1, which
establishes then the Abelian result, since the sum
in (2.1) is seen to be asymptotic to (2.2), similarly as was done
in [6], and thus by a change of variable it follows that (2.3) holds.
We have
$$
\eqalign{
g(\lambda x) - g(x) &= \int_1^{\lambda A}{[u]\over u^2}
{\tilde \ell}\left({\lambda x\over u}\right)\d u
- \int_1^A{[u]\over u^2}
{\tilde \ell}\left({x\over u}\right)\d u \cr&
+ \int_A^{x/2}\left({\{u\}\over u^2} - {\{\lambda u\}\over \lambda u^2}\right)
{\tilde \ell}\left({x\over u}\right)\d u,\cr}
$$
where $\{u\} = u - [u]$ is the fractional part of $u$ and $A$
is a large positive constant. By the uniform convergence theorem we obtain
$$
\eqalign{&
\int_1^{\lambda A}{[u]\over u^2}
{\tilde \ell}\left({\lambda x\over u}\right)\d u - \int_1^A{[u]\over u^2}
{\tilde \ell}\left({x\over u}\right)\d u
= (1+o(1)){\tilde\ell}(x)\int_A^{\l A}
{[u]\over u^2}\d u\cr&
= \left(1+o(1)+O\left({1\over A}\right)\right){\tilde\ell}(x)
\int_A^{\l A}{\d u\over u}
= \left(\log\l + o(1)+ O\left({1\over A}\right)\right){\tilde\ell}(x).
\cr}
$$
We also have
$$
\eqalign{&
\int_A^{x/2}\left({\{u\}\over u^2} - {\{\lambda u\}\over \lambda u^2}\right)
{\tilde \ell}\left({x\over u}\right)\d u
\cr&
= (1+o(1)){\tilde\ell}(x)\left(\int_A^\infty {\{u\}\over u^2}\d u
- \int_A^\infty {\{\lambda u\}\over \lambda u^2}\d u\right) \ll
{{\tilde\ell}(x)\over A},
\cr}
$$
hence by letting $A\to\infty$ it follows that $g\in\Pi_{\tilde\ell}$
with ${\tilde\ell}$--index 1. Here we used the fact that, for $L\in{\Cal L}$,
we have
$$
\int_{x_0}^{Bx} L\left({x\over u}\right)h(u)\d u \;\sim\;
L(x)\int_{x_0}^\infty h(u)\d u\qquad(x\to \infty)\leqno(2.6)
$$
if $B>0$ is a constant, and $h$ is an integrable function
satisfying $h(u) \ll u^{-c}$ for some $c>1$. One obtains (2.6) by writing
$$
\int_{x_0}^{Bx} = \int_{x_0}^A +
\int_A^{x/A} + \int_{x/A}^{Bx} = I_1 + I_2 + I_3,
$$
say, where again $A$ is a large constant. By the uniform convergence theorem
$$
I_1 = (1+o(1))L(x)\int_{x_0}^A h(u)\d u =
\left(\int_{x_0}^\infty h(u)\d u+o(1) + O\left(A^{1-c}\right)\right)L(x).
$$
Next, since $x^{-\e}L(x)$ (see (1.6)) is asymptotic to a non-increasing
function for $x \geqslant x_1(\e)$ and any given $\e>0$, we have
$$\eqalign{
I_2 &= \int_A^{x/A} \left({x\over u}\right)^{-\e}L\left({x\over u}\right)
\left({x\over u}\right)^{\e}
h(u)\d u \cr&\ll \left({x\over A}\right)^{-\e}L\left({x\over A}\right)
\int_A^{\infty}\left({x\over u}\right)^{\e}|h(u)|\d u \ll
L\left({x\over A}\right)A^{1-c} = (A^{1-c} + o(1))L(x).\cr}
$$
Finally
$$
I_3 \ll \int_{x/A}^\infty u^{-c}\d u \;\ll\; {\left({x\over A}\right)}^{1-c},
$$
and (2.6) follows on taking $A$ sufficiently large, since $c>1$.

\medskip
The above analysis shows that the case $\rho = 0$ is more delicate than
the case $\rho >0$. Namely the function $g(x)$ given by (2.3)--(2.4)
is not so simple to evaluate asymptotically.
To illustrate our point, consider  the additive function
$$
G_\a(n) := \sum_{p|n}\log^\a p\qquad(\a \in\RR),
$$
which is of the form (1.1) with $\rho=0$ and $L(x) = \log^\a x$
($G_0(n) \equiv \omega(n)$). Then
$$
\sum_{n\leqslant x} G_\a(n) =  \sum_{pm\leqslant x}\log^\a p
= \sum_{p\leqslant x}\left[{x\over p}\right]\log^\a p.     \leqno(2.7)
$$
We use the prime number theorem in its strongest known form, namely
$$
\pi(x) :=  \sum_{p\leqslant x}1 = \int_2^x{\d t\over\log t} + R(x)
$$
with
$$
R(x) \ll x{\roman e}^{-c\e(x)},\quad \e(x) := \log^{3/5}x(\log\log x)^{-1/5},
\quad c > 0.
$$
Then we have
$$ \eqalign{&
\sum_{p\leqslant x}\left[{x\over p}\right]\log^\a p =
\int_{2-0}^x \left[{x\over t}\right](\log t)^{\a}\d \pi(t)\cr&
= \int_2^x \left[{x\over t}\right](\log t)^{\a-1}\d t +
\int_{2-0}^x \left[{x\over t}\right](\log t)^{\a}\d R(t).\cr}\leqno(2.8)
$$
One can show, similarly as was done in [6], that
$$
\int_{2-0}^x \left[{x\over t}\right](\log t)^{\a}\d R(t)
\ll x{\roman e}^{-{1\over2}c\e(x)}.
$$
We have
$$
\int_2^x \left[{x\over t}\right](\log t)^{\a-1}\d t
=   x\int_2^x  t^{-1}(\log t)^{\a-1}\d t + O(x(\log x)^{\a-1}) \leqno(2.9)
$$
and
$$
\int_2^x  t^{-1}(\log t)^{\a-1}\d t
\,=\, \Bigg\{\aligned {\log^\a x\over\a} + O(1)&\quad(\a > 0), \\
\log\log x + O(1)&\quad(\a = 0),
\\
\int_2^\infty (\log t)^{\a-1}{\d t\over t} + &O((\log x)^{\a})
\;(\a < 0).\endaligned
\leqno(2.10)
$$
The asymptotic formula for the summatory function of $G_\a(n)$ follows
then from (2.7)-(2.10). Its shape changes according to the cases
$\a < 0,\,\a=0$ or $\a > 0$, respectively.

\medskip
It was pointed out in [6], without details of proof, that (2.1) remains
valid if $f(n)$ is replaced by $F(n)$ (both given by (1.1)). To see this,
note first that for any additive function $g(n)$ one has
$$
\eqalign{&
\sum_{n\leqslant x}g(n) = \sum_{n\leqslant x}\sum_{p^\nu||n}g(p^\nu)
= \sum_{p^\nu m\leqslant x,(p,m)=1} g(p^\nu)\cr&
= \sum_{p^\nu m\leqslant x} g(p^\nu)
- \sum_{p^{\nu+1} m\leqslant x} g(p^\nu) =    \sum_{p^\nu \leqslant x}
(g(p^\nu) - g(p^{\nu-1}))\left[{x\over p^\nu}\right]\cr&
=   \sum_{p\leqslant x}g(p)\left[{x\over p}\right]
+   \sum_{p^\nu \leqslant x,\nu\geqslant2}
(g(p^\nu) - g(p^{\nu-1}))\left[{x\over p^\nu}\right].
\cr} \leqno(2.11)
$$
Setting $g(n) = F(n) - f(n)$ in (2.11) and noting that
$$
g(p) = 0,\quad g(p^\nu) = (\nu-1)p^\rho L(p) \;(\nu \geqslant 2),
$$
it follows that
$$
\eqalign{&
\sum_{n \leqslant  x}(F(n) - f(n)) =
\sum_{p^\nu \leqslant x,\nu\geqslant2} p^\rho L(p)
\left[{x\over p^\nu}\right]\cr&
\leqslant  x\sum_{p\leqslant \sqrt{x}}p^\rho L(p)\left(
{1\over p^2} +  {1\over p^3} + \ldots\,\right) \ll
x\sum_{p\leqslant \sqrt{x}}p^{\rho-2}L(p).\cr}\leqno(2.12)
$$
Since $x^{-\e} \ll L(x) \ll x^\e$ for any given $\e>0$, we have
$$
\sum_{p\leqslant \sqrt{x}}p^{\rho-2}L(p)
\ll \sum_{p\leqslant \sqrt{x}}p^{\rho-2+\e} \ll_\e
\Bigg\{\aligned 1&\quad(\rho < 1),
\\    x^{{1\over2}(\rho-1)+\e}
&\quad(\rho \geqslant 1).\endaligned\leqno(2.13)
$$
Hence for sufficiently small $\e$ we obtain, in view of (2.1),
$$
\eqalign{
\sum_{n \leqslant  x}F(n) &= \sum_{n \leqslant  x}f(n)
+ O_\e\left(\max(x,\,x^{{1\over2}(\rho+1)+\e}) \right)\cr&
= \left({\zeta(\rho+1)\over\rho+1} + o(1)\right){x^{1+\rho}L(x)\over\log x}
\quad(\rho>0,\,x\to\infty),\cr}
$$
which proves our assertion. In the case $\rho = 0$ the asymptotic formula
(2.3)-(2.4) remains valid if and only if
$$
\lim_{x\to\infty}g(x) = +\infty.\leqno(2.14)
$$
To see this note that if (2.14) holds, then (2.12) gives
(in case $\rho=0$) that
$$
\sum_{n \leqslant  x}(F(n) - f(n)) \ll x,
$$
hence
$$
\sum_{n \leqslant  x}F(n) = (1+o(1))xg(x)\qquad(x\to\infty).\leqno(2.15)
$$
On the other hand, if $g(x)$ is bounded, then (2.15) cannot hold.
Namely for $\rho=0$ we obtain from (2.12)
$$
\sum_{n \leqslant  x}F(n) = \sum_{n \leqslant  x}f(n) +
\sum_{p^\nu \leqslant  x,\nu\geqslant 2}L(p)\left[{x\over p^\nu}\right].
$$
Then the last sum above equals
$$
x\sum_{p^\nu \leqslant  x,\nu\geqslant 2}L(p)p^{-\nu} +
O_\e(x^{{1\over2}+\e}) = Cx + O_\e(x^{{1\over2}+\e})
$$
with
$$
C \;:=\; \sum_p {L(p)\over p(p-1)} \;(< +\infty).
$$
This means that
$$
\sum_{n \leqslant  x}F(n) = x(g(x) + C + o(1))\qquad(x\to\infty),
$$
which proves our claim.

\heading
3. Tauberian theorems for sums of additive functions
\endheading

Tauberian theorems for arithmetic sums, including the Tauberian
counterpart of (2.1) (both the case $\rho > 0$ and the limiting case
$\rho=0$), have been recently obtained by Bingham--Inoue [3]. In
previous works (see e.g., [2]) they developed powerful Tauberian theorems
for systems of kernels. In this way they succeeded in simplifying
and extending Tauberian results of De Koninck and the author [8].
In the case  $\rho > 0$ they proved that
$$
\sum_{n\leqslant x}\sum_{p|n}h(p) = \left({\z(\rho+1)\over\rho+1}
+ o(1)\right){x^{\rho+1}L(x)\over\log x}\quad(x\to\infty)\leqno(3.1)
$$
with $L\in{\Cal L}$, $h$ positive and non-decreasing  implies that
$$
h(x) \;\sim\; x^\rho L(x)\qquad(x\to\infty).\leqno(3.2)
$$
This result they call Tauberian, and more generally they proved that if
$$
\sum_{n\leqslant x}\sum_{p|n}h(p) = (C + o(1))h(x){x\over\log x}
\quad(x\to\infty)\leqno(3.3)
$$
holds for some positive constant $C$, and $h\;:\;[2,\infty)\;\rightarrow
\;(0,\infty)\,$ is continuous and non-decreasing, then $h$ is regularly
varying with index $\rho$, where $\rho$ is the unique positive solution
of the equation $C = \z(\rho+1)/(\rho + 1)$. This type of result they called
{\it Mercerian}. In the limiting case $\rho=0$ they prove the
following Tauberian theorem: If $L\in{\Cal L}$,
$$
\int_{x_0}^\infty L(x){\roman e}^{-\sqrt{\log x}}{\d x\over x }
\;<\;\infty,\leqno(3.4)
$$
$$
\log x \;=\;O(L(x)),\leqno(3.5)
$$
then if ${\tilde L}(x) = L(x)/\log x$ and
$$
{1\over x}\sum_{n\leqslant x}\sum_{p|n}h(p) \;\in\;\Pi_{\tilde L}
\qquad(\text{with}\,{\tilde L}-\text{index}\; 1)\leqno(3.6)
$$
with $h\;:\;[2,\infty)\;\rightarrow\;(0,\infty)\,$  continuous and
non-decreasing, we have
$$
h(x) \;\sim\; L(x)\qquad(x\to\infty).\leqno(3.7)
$$
The Tauberian condition (3.4) comes from using the prime theorem
with the classical error term (see Section 2) $R(x) \ll
x{\roman e}^{-\sqrt{\log x}}$. It could be replaced with a slightly
sharper condition coming from the strongest known error term, namely
$$
\int_{x_0}^\infty L(x){\roman e}^{-c(\log x)^{3/5}(\log\log x)^{-1/5}}
{\d x\over x } \;<\;\infty,
$$
with suitable $c > 0$.
However, as remarked by Bingham-Inoue [3], the Tauberian condition
(3.4) is much less restrictive than (3.5). It is an open problem to
relax (3.5), or dispense with it altogether. The method of proof of the
above results consists of using Tauberian (Mercerian) theorems for
the {\it Mellin convolution}
$$
(f\star k)(x) := \int_0^\infty k\left({x\over t}\right)f(t){\d t\over t},
$$
on showing that the relevant sums are asymptotic to an appropriate Mellin
convolution.

\medskip
One can also ask for the Tauberian analogue of (3.1) if we suppose
that $$ \sum_{n\leqslant x}\sum_{p^\nu||n}\nu h(p) =
\left({\z(\rho+1)\over\rho+1} +
o(1)\right){x^{\rho+1}L(x)\over\log x}\quad(x\to\infty)\leqno(3.8)
$$ with $\rho>0$, $h$ positive and non-decreasing. Note that $$
\eqalign{& \sum_{n\leqslant 2x}\sum_{p^\nu||n}\nu h(p) \geqslant
\sum_{n\leqslant 2x}\sum_{p|n}h(p) = \sum_{p\leqslant
2x}h(p)\left[{ 2x\over p}\right] \cr& \geqslant \sum_{x<p\leqslant
2x}h(p) \geqslant h(x)(\pi(2x)- \pi(x)) \gg {xh(x)\over\log
x},\cr}\leqno(3.9) $$ and consequently (3.8) yields $$ h(x)
\;\ll\; x^\rho L(x).\leqno(3.10) $$ Therefore by using (2.12),
(2.13) and (3.10) we find that $$ \eqalign{ \sum_{n\leqslant
x}\sum_{p^\nu||n}\nu h(p)  &= \sum_{n\leqslant x}\sum_{p|n}h(p) +
\sum_{p^\nu\leqslant x,\nu\geqslant2} h(p)\left[{x\over
p^\nu}\right] \cr& = \sum_{n\leqslant x}\sum_{p|n}h(p) + O\Bigl(
\sum_{p^\nu\leqslant x,\nu\geqslant2}p^\rho L(p)\left[{x\over
p^\nu}\right] \Bigr)\cr& = \sum_{n\leqslant x}\sum_{p|n}h(p) +
O\Bigl(x\sum_{p\leqslant\sqrt{x}}p^{\rho-2} L(p)\Bigr)\cr& =
\sum_{n\leqslant x}\sum_{p|n}h(p) +
O_\e\left(\max(x,\,x^{{1\over2}(\rho+1)+\e})\right).
\cr}\leqno(3.11) $$ Since $\rho > 0$ it follows that (3.8) yields
(3.1), hence by the Tauberian theorem of Bingham--Inoue we have
(3.2). A similar discussion shows that the corresponding Mercerian
analogue of (3.3) also holds for the sum in (3.8), as well as does
the Tauberian anlogue of (3.6). For the latter note that, by the
analogue of (3.9) when $\rho = 0$, we obtain that $$ h(x) \ll
g(x)\log x,\quad g(x) \in \Pi_{\tilde L}
\qquad(\text{with}\,{\tilde L}-\text{index}\; 1), $$ and
consequently $h(x) \ll_\e x^\e$. Supposing $\l > 1$ (the case $\l
< 1$ is analogous), we have $$ \eqalign{& {1\over\l
x}\sum_{p^\nu\leqslant \l x,\nu\geqslant 2}h(p) \left[{\l x\over
p^\nu}\right] - {1\over x}\sum_{p^\nu\leqslant x,\nu\geqslant
2}h(p) \left[{x\over p^\nu}\right] \cr& = \sum_{x<p^\nu\leqslant
\l x,\nu\geqslant 2}{h(p)\over p^\nu} + O\left({1\over
x}\sum_{p^\nu\leqslant \l x,\nu\geqslant 2} {|h(p)|\over
p^\nu}\right) = O_\e(x^{\e-{1\over2}}). \cr} $$ In view of (3.6)
and (3.11) this means that $$ \sum_{n\leqslant
x}\sum_{p^\nu||n}\nu h(p) \in \Pi_{\tilde L}
\qquad\left(\text{with}\,{\tilde L}-\text{index}\; 1,\,{\tilde
L}(x) = {L(x)\over\log x}\right),
$$ implies  (3.6), and therefore
(3.7) follows if (3.4)-(3.5) holds and  the function
$h\;:\;[2,\infty)\;\rightarrow\;(0,\infty)\,$ is continuous and
non-decreasing.

\heading
4. Sums of reciprocals
\endheading
We shall conclude our exposition by considering sums of reciprocals
of certain arithmetic functions.
It is a classical result of prime number theory that
$$
\sum_{p\leqslant x}{1\over p} = \log\log x + C + O\left({1\over\log x}\right).
$$
Sums of reciprocals of large additive functions $\b(n),\,B(n)$ (see (1.8))
and $P(n)$  are much
more difficult to handle, where $P(n)$ is the largest prime factor
of $n \,(\geqslant 2)$ and $P(1) = 1$. They were
investigated by De Koninck, Erd\H os, Pomerance, Xuan and the author
(see [4]--[7], [9]--[14], [16]--[19], [22], [23]).
It was proved by  Erd\H os, Pomerance and the author [13] that
$$
\sum_{n\leqslant x}{1\over P(n)} = x\delta(x)\left(1 + O\left(\sqrt{{\log\log x
\over\log x}}\,\right)\right)\leqno(4.1)
$$
with
$$ \delta(x) \;:=\; \int_2^x \rho\left(
{\log x\over\log t}\right){\d t\over t^2},\leqno(4.2)
$$
where the Dickman--de Bruijn function $\rho(u)$
is the continuous solution to the differential delay equation
$$u\rho'(u) = -\rho(u-1), \;\rho(u) = 1\;
\text{for}\; 0 \leqslant u \leqslant 1,\; \rho(u) = 0\;\text{for}\; u < 0\;.
$$
It is known (see [21]) that ($\log_kx = \log(\log_{k-1}x$))
$$
\rho(u) = \exp\Biggl\{-u\Biggl(\log u + \log_2u - 1
+ {\log_2u-1\over\log u} + O\left(\left({\log_2u\over\log u}\right)^2
\right)\Biggr)\Biggr\}.
$$
\medskip
It was shown by Pomerance and the author [19] that one has
$$
\delta(x) = \exp\left\{-(2\log x\log_2x)^{1/2}
\left(1 + g_0(x) + O\left({\log_3^3x\over\log_2^3x}\right)\right)\right\},
$$
where
$$\eqalign{
g_r(x) &= {\log_3x + \log(1+r) -2 -\log2\over2\log_2x}\left(1 +
{2\over\log_2x}\right) \cr&
- \,{(\log_3x + \log(1+r) - \log2)^2\over8\log_2^2x},\cr}
$$
and the expression for $\delta(x)$ was
sharpened by the author [17]. Already in 1977 Paul Erd\H os told the
author that the function $\delta(x)$ is slowly varying,
but it is only in 1986 that this fact was established in [13], by the
use of (4.2) and properties of the function $\rho(u)$. A corollary of
(4.1) and the fact that $\delta(x)$ is slowly varying is the asymptotic
formula
$$
\sum_{n\leqslant x}{1\over P(n)} \;\sim\;\sum_{x<n\leqslant 2x}{1\over P(n)}
\qquad(x\to\infty),
$$
which is by no means obvious. In fact we have (see [13, p. 291]) that
$$
\delta\left({x\over t}\right) = \left(1 + O\left(
{(\log\log x)^{3/2}\over(\log x)^{1/2}}\right)\right)\delta(x)
\qquad(1 \leqslant t \leqslant \log^{12}x).
$$
One may ask: precisely for what $t\geqslant 1$ does one have
$\delta(x/t) \sim \delta(x)$ as $x\to\infty$?

\medskip
The asymptotic formula (4.1) remains valid if $P(n)$ is replaced by
$\b(n)$ or $B(n)$, and the asymptotic formula for the summatory
function of ${B(n)\over\b(n)}-1$ is of the same shape as the right-hand
side of (4.1). Furthermore we have
$$\eqalign{&
\sum_{2\leqslant n\leqslant x}\left({1\over\b(n)} - {1\over B(n)}\right) \cr&
=  x\exp\left\{-2(\log x\log_2x)^{1/2}\left(1 + g_1(x)
+ O\left({\log_3^3x\over\log_2^3x}\right)\right)\right\}.\cr}\leqno(4.3)
$$
One can also show that  the sum on the left-hand side of (4.3)
is asymptotic to a regularly varying function with index 1.

\medskip
Finally let us mention that,
by using results of the joint paper with   Erd\H os and  Pomerance [13],
the author [16] sharpened some of the  asymptotic formulas proved in
earlier works and obtained (see (1.7)), for example,
$$
\sum_{n\leqslant x}{\Omega(n)-\omega(n)\over P(n)}
= \left\{\sum_p{1\over p^2-p}+ O\left(\left({\log_2x\over\log x}\right)
^{1/2}\right)\right\}\sum_{n\leqslant x}{1\over P(n)},
$$
$$
\sum_{n\leqslant x}{\omega(n)\over P(n)} =
 \left\{\left({2\log x\over\log_2x}\right)^{1/2}
\left(1 +  O\left({\log_3x\over\log_2 x}\right)
\right)\right\}\sum_{n\leqslant x}{1\over P(n)},
$$
and this remains valid if $\omega(n)$ is replaced by $\Omega(n)$.
We also have
$$
\sum_{n\leqslant x}{\mu^2(n)\over P(n)} =
\left\{{6\over\pi^2}+O\left(\left({\log_2x\over\log x}\right)^{1/2}
\right)\right\}\sum_{n\leqslant x}{1\over P(n)},  \leqno(4.4)
$$
where $\mu$ is the M\"obius function.  On the left hand-side
of (4.4) we have summation over squarefree numbers, whose density
is precisely $6\over\pi^2$, so that this result shows that the
sum of reciprocals of $P(n)$ over squarefree numbers $\leqslant x$
asymptotically equals the corresponding density times the sum of
$1/P(n)$ over all $n\leqslant x$.

\bigskip\bigskip

\Refs

\item{[1]} N.H. Bingham, C.M. Goldie and J.L. Teugels,
Regular Variation, Cambridge University Press, Cambridge, 1987.

\item{[2]} N.H. Bingham and A. Inoue, Tauberian and Mercerian theorems
for systems of kernels, J. Math. Anal. Appl. {\bf252}(2000), 177-197.

\item{[3]} N.H. Bingham and A. Inoue,   Abelian,
Tauberian and Mercerian theorems for arithmetic sums,
J. Math. Anal. Appl. {\bf250}(2000), 465-493.

\item{[4]} J.-M. De Koninck and A. Ivi\'c, Topics in arithmetical
    functions, Mathematics Studies {\bf43}, North-Holland,
    Amsterdam, 1980.

\item{[5]} J.-M. De Koninck and A. Ivi\'c,
 The distribution of the average prime
divisor of an integer, Archiv  Math.
{\bf43}(1984), 37-43.

\item{[6]} J.-M. De Koninck and A. Ivi\'c, On the average prime factor
of an integer and some related problems, Ricerche di Matematica
{\bf39}(1990), 131-140.

 \item{[7]} J.-M. De Koninck and A. Ivi\'c, Random sums related to prime
 divisors of an integer, Publs. Inst. Math. (Beograd) {\bf48(62)}(1990),
 7-14.

 \item{[8]} J.-M. De Koninck and A. Ivi\'c, Arithmetic characterization
 of regularly varying functions,  Ricerche di Matematica {\bf44}(1995),
 41-64.

\item{[9]} J.-M. De Koninck, P. Erd\H os and A. Ivi\'c,
    Reciprocals of certain large additive
    functions, Canadian Math. Bulletin {\bf24}(1981), 225-231.

\item{[10]} P. Erd\H os and A. Ivi\'c,
Estimates for sums involving the largest prime     factor of an
    integer and certain related additive functions, Studia Scientiarum
    Math. Hungarica {\bf15}(1980), 183-199.

\item{[11]}  P. Erd\H os and A. Ivi\'c,
On sums involving reciprocals of certain arithmetical functions,
Publications Inst. Math. (Beograd) {\bf 32(46)}(1982), 49-56.

\item{[12]}  P. Erd\H os and A. Ivi\'c,
 The distribution of small and large additive functions II,
 Proceedings of the Amalfi Conference on Analytic Number Theory
    (Amalfi, Sep. 1989), Universit\`a di Salerno, Salerno 1992, 83-93.

\item{[13]} P. Erd\H os, A. Ivi\'c and C. Pomerance,
On sums involving reciprocals of the
largest prime factor of an integer,  Glasnik Matemati\v cki
{\bf21(41)}(1986), 283-300.

\item{[14]} A. Ivi\'c, Sum of reciprocals of the largest
prime factor of an integer, Archiv
    Math. {\bf36}(1981), 57-61.

\item{[15]} A. Ivi\'c, The Riemann zeta-function, John Wiley \& Sons,
New York, 1985.

\item{[16]} A. Ivi\'c, On some estimates involving the
number of prime divisors of an integer,
    Acta Arithmetica {\bf49}(1987), 21-32.

\item{[17]} A. Ivi\'c, On sums involving reciprocals of
the largest prime factor of an integer II,
        Acta Arithmetica {\bf71}(1995), 241-245.

\item{[18]} A. Ivi\'c, On large additive functions over
primes of positive density,
        Mathematica Balkanica {\bf10}(1996), 103-120.

\item{[19]} A. Ivi\'c and C. Pomerance,
 Estimates for certain sums involving the
    largest prime factor of an integer, Proceedings Budapest Conference
    in Number Theory July 1981, Coll. Math. Soc. J. Bolyai {\bf34},
    North-Holland, Amsterdam 1984, 769-789.

\item{[20]} J. Karamata, Sur un mode de croissance r\'eguli\`ere des
fonctions, Mathematica (Cluj) {\bf4}(1930), 38-53.

\item{[21]} E. Seneta, Regularly varying functions, LNM {\bf508}, Springer
Verlag, Berlin--Heidelberg--New York, 1976.

\item{[21]} G. Tenenbaum, Introduction \`a la th\'eorie analytique
et probabiliste des nombres, Soci\'et\'e Math. de France, Paris, 1995.

\item{[22]} T.Z. Xuan, On sums involving reciprocals of certain large
additive functions, Publs. Inst. Math. (Beograd) {\bf45(59)}, 41-55
and II, ibid. {\bf46(60)}(1989), 25-32.

\item{[23]} T.Z. Xuan, On a result of Erd\H os and Ivi\'c, Archiv Math.
{\bf62}(1994), 143-154.

\bigskip

 Aleksandar Ivi\'c

Katedra Matematike RGF-a

Universitet u Beogradu

\DJ u\v sina 7, 11000 Beograd, Serbia

aivic\@rgf.bg.ac.yu, aivic\@matf.bg.ac.yu

\endRefs

%\enddocument
\bye